\documentclass[11pt]{amsart}
\usepackage{amsmath,amssymb,mathrsfs,stmaryrd}
\newtheorem{theorem}{Theorem}

\newtheorem{example}[theorem]{Example}
\newcommand{\CP}{{\mathbb C}{\mathbb P}}
\newcommand{\R}{\mathbb R}
\newcommand{\C}{\mathbb C}

\begin{document} 

\title[Solitons solutions of geometric flows]{Constructing soliton solutions of geometric flows by separation of variables}
\author{Mu-Tao Wang}
\begin{abstract}
This note surveys and compares results in \cite{Futaki-Wang} and \cite{Lee-Wang1, Lee-Wang2} on the separation of variables construction for soliton solutions of curvature equations including the K\"ahler-Ricci flow and the Lagrangian mean curvature flow. In the last section, we propose some new generalizations in the Lagrangian
mean curvature flow case.  
\end{abstract}
\address{Department of Mathematics \\ Columbia University \\ 2990 Broadway \\ New York, NY 10027}
\thanks{The author was supported in part by the National Science Foundation under grant DMS-1105483.}
\maketitle 

\section{Introduction}

 Let us recall a simple example of separation of variables constructions of solutions of elliptic PDE's. Consider the Laplace equation $\Delta^{\R^n} f=0$ on $\R^n$. In terms of polar coordinates, the equation becomes
\begin{equation}\label{polar}\frac{1}{r}\frac{\partial}{\partial r}(r^2\frac{\partial f}{\partial r})+\Delta^{S^{n-1}} f=0.\end{equation}

The separation of variables ansatz assumes that $f$ takes the form $f (r, \theta) =f_1(r) f_2(\theta)$ for $r\in (0, \infty)$ and $\theta \in S^{n-1}$, and  equation \eqref{polar}
can be written as 
\[-\frac{\frac{1}{r}\frac{\partial}{\partial r}(r^2\frac{\partial f_1}{\partial r})}{f_1(r)}=\frac{\Delta^{S^{n-1}} f_2}{f_2}=\text{constant}.\]
Thus, the original equation reduces to an ODE for $f_1$ and a PDE for a function $f_2$  on the sphere.

Suppose instead we would like  to solve an inhomogeneous equation of the form $\Delta^{\R^n} f=g$. Assuming again that $g$ takes the form $g(r,\theta)=g_1(r) g_2(\theta)$, for $g_1 (r)$ a function on $(0, \infty)$ and $g_2(\theta)$ a function on $S^{n-1}$, the ansatz $f(r,\theta)=f_1(r) f_2(\theta)$ for suitable choices of $g_2(\theta)$ and $f_2(\theta)$ reduces the original equation to an ODE for $f_1(r)$ involving $g_1(r)$. Thus the method of separation of variables produces solutions of homogeneous equations and reduces an inhomogeneous equation to an ODE.

\section{Application to nonlinear elliptic curvature equations}
We discuss two types of geometric parabolic equations, the Lagrangian mean curvature flow and the K\"ahler-Ricci flow, and their soliton solutions. First we review the elliptic case. 

Suppose a submanifold is given by an embedding $X$ (or more generally an immersion) into the Euclidean space. The mean curvature equation $H(X)=0$ corresponds to minimal submanifolds. We recall that $H(X)=\Delta X$ for the embedding $X$ where $\Delta$ is the Laplace operator with respect to the induced metric. Thus this should be considered as the Laplace equation for submanifolds.

On the other hand, suppose a Riemannian metric $g$ is given.  The Ricci curvature equation $Ric (g)=0$ corresponds to Ricci flat manifolds. We recall that $Ric(g)= -\frac{1}{2} \Delta g$ at a point with respect to a harmonic normal coordinate system, and we can view the Ricci equation as a Laplace equation for Riemannian metrics.

For both nonlinear elliptic equations, cone-like solutions appear naturally in the study of singularity formation through blow-up, or rescaling process. Minimal cone in $\R^n$ is of the form $rX, r>0$ where $X:\Sigma \rightarrow S^{n-1}$ is a minimal submanifold in $S^{n-1}$. On the other hand, a Ricci flat cone is of the form $dr^2+r^2ds^2_\Sigma$ for an Einstein manifold  $\Sigma$. Both are solutions of homogeneous elliptic equations and provide ansatz for separation of variable. An elliptic geometric curvature problem is often accompanied by a parabolic problem,  which is the gradient flow in case the elliptic problem is variational. Parabolic blow-up (rescaling) gives soliton solutions  which are self-similar in space and time variables. To classify singularities of parabolic flows,  it is important to study solitons. It turns out solitons of parabolic equations satisfy elliptic inhomogeneous equations.

In the case of mean curvature flows, 
suppose $X_t$ is a family of embeddings that satisfies the mean curvature flow equation: \[\frac{\partial}{\partial t} X_t=H (X_t).\] Special solutions called  solitons or  self-similar solutions can be obtained in the following way. First one solves the elliptic equation for an embedding $X$ such that \[H(X)=(aX+\vec{b})^\perp\] for a constant real number $a$ and a constant vector $\vec{b}$.  Then one defines a one-parameter family of embeddings $X_t$ by
\[X_t=a(t)X+\vec{b}(t).\] One can check that $X_t$ forms a mean curvature flow that is moving by scaling and translation for suitable choices of $a(t)$ and $\vec{b}(t)$. The sign of $a$ determines whether the flow is expanding, shrinking, or steady (translating).

 On the other hand, the Ricci flow is a parabolic equation for a family of metrics $g_t$ that satisfies: \[\frac{\partial}{\partial t} g_t=-2 Ric (g_t).\] Self-similar solutions can be obtained by solving the elliptic  Ricci soliton equation: $Ric(g) =a g+ \mathcal{L}_V g$ for a metric $g$, where $a$ is $a$ constant real number and $V$ is a smooth vector field on the underlying manifold.

The family of metrics $g_t$
 \[g_t=a(t)\gamma_t^* g,\] where $\gamma_t$ is the one-parameter group of diffeomorphism generated by $V$,  satisfies the Ricci flow equation and the corresponding flow is a combination of scalings and reparametrizations. Likiwise, the sign of $a$ determines whether the flow is expanding, shrinking, or steady.

\section{The K\"ahler-Ricci flow}

Let us first look at the K\"ahler-Ricci flow. Suppose the K\"ahler form is $\omega=i g_{i{\bar{j}} }dz^i\wedge d\bar{z}^j$, the Ricci form is then given by $\rho(\omega)=-i\partial\bar{\partial}\log \det g_{i\bar{j}}$. Given a family of K\"ahler metrics, in terms of the associated K\"ahler forms, the K\"ahler-Ricci flow is the parabolic equation \[\frac{d}{dt}\omega=-\frac{1}{2}\rho(\omega).\] Correspondingly,  the K\"ahler-Ricci soliton equation is 
\[-\frac{1}{2}{\rho}(\omega)=\lambda\omega+\mathcal{L}_V \omega,\]  where $\lambda=1, 0, -1$ corresponds to expanding, steady and shrinking Ricci solitons.

This is a gradient soliton if \[\mathcal{L}_V \omega=i\partial \bar{\partial} Q\] for a function $Q$. Therefore, the  K\"ahler-Ricci gradient soliton equation for K\"ahler form $\omega$ and potential $Q$ is:
\[ -\frac{1}{2}\rho(\omega)=\lambda\omega+i\partial\bar{\partial} Q.\]

The Ansatz to apply the separation of variables method in this case involves a Sasakian manifold. Recall a $2m+1$ dimensional Riemannian manifold $(S, g)$ is Sasakian if $\bar{g}=dr^2+r^2g$ is a K\"ahler metric on the cone $C(S)$ over $S$.

The restriction to $S$ of ${\xi}=Jr\frac{\partial}{\partial r}$ is the Reeb vector field  of the contact form that satisfies $i(\xi)\eta=1$ and $i(\xi)d\eta=0$. We have the following relation:  \[\eta=\frac{1}{r^2} \bar{g}(\xi, \cdot)=i(\bar{\partial}-\partial)\log r.\]

The vector field $\xi$ defines a transverse K\"ahler structure on the local orbit space.  $\xi$ is Killing on $(C(S), \bar{g})$ and its complexification $\xi-iJ\xi$ is holomorphic.   In particular, the transverse metric satisfies $g^T+\eta\otimes \eta=g$ and the transverse K\"ahler form is
\[\omega^T=\frac{1}{2} d\eta=i\partial\bar{\partial} \log r.\] The K\"ahler form on $C(S)$ is \[\omega=\frac{i}{2}\partial \bar{\partial} r^2.\]
For example, we can take $S=S^{2m+1}$, $C(S)=\C^{m+1}-\{0\}$ and the orbit space is $\CP^{m}$.

One can check that $(S, g)$ is a Sasaki-Einstein manifold if and only if $(C(S), \bar{g})$ is a Calabi-Yau cone manifold. Also, $S$ is $\eta$-Einstein ($Ric_g=\alpha g+\beta \eta\otimes \eta$ with $\alpha+\beta=2m$) if and only if the transverse K\"ahler metric is Einstein $(\rho(\omega^T)=(\alpha+2)\omega^T)$.

Suppose we start with an $\eta$-Einstein Sasaki manifold. The ansatz for the K\"ahler form of a K\"ahler-Ricci gradient soliton is
\[\omega^T+i\partial\bar{\partial}F(s)\] with $s=\log r$. This is similar to the Calabi Ansatz \cite{Calabi} for Calabi-Yau metrics which was applied to the total space of a Hermitian line bundles over a K\"ahler manifold.

The equation for $F(s)$ can be derived as follows.  Set $\sigma=1+ F'(s)$ and $\phi(\sigma)=F''(s)$. We have
\[\omega_\phi=\sigma\omega^T+\phi(\sigma) i \partial s\wedge \bar{\partial}s\]
where $i \partial s\wedge \bar{\partial}s$ is the K\"ahler form of the cylindrical metric on $\C-\{0\}$.

As a comparison, the metric on the Calabi-Yau cone ($Ric=\rho(\omega)=0$) is given by $F(s)=\frac{r^2}{2}-\log r=\frac{1}{2} e^{2s}-s$ and $\phi(\sigma)=2\sigma$. The Ricci form  of $\omega_\phi$ is\[\rho (\omega_\phi)=\kappa \omega^T-i\partial\bar{\partial}\log(\sigma^m\phi(\sigma))\] with $\kappa=\alpha+2$.

For a gradient soliton, if the potential $Q$ is a function of $s$, $Q$ must be of the form $Q=\mu\sigma+c$ for constants $\mu$ and $c$. Therefore, the gradient soliton equation becomes
\[\phi'(\sigma)+(\frac{m}{\sigma}-\mu)\phi(\sigma)-(\kappa+2\lambda \sigma)=0\] for $m, \mu, \kappa,\lambda$ constants.

Studying of the ODE gives the following theorem:
\begin{theorem} \cite{Futaki-Wang} Let $S$ be a compact Sasaki manifold such that the transverse K\"ahler metric $g^T$ satisfies
Einstein equation
$$ \mathrm{Ric}^T = \kappa g^T$$
for some $\kappa < 0$ (or $S$ is $\eta$-Einstein). Then there exists a complete
expanding soliton on the K\"ahler cone $C(S)$.
 \end{theorem}

On total space of line bundles over Fano manifolds, we also obtain shrinking and expanding solitons that can be glued together to form eternal solution on $(-\infty, \infty)$.
\begin{theorem}\cite{Futaki-Wang} 
 Let $M$ be a Fano manifold of dimension $m$, and $L \to M$ be a positive line bundle with $K_M = L^{-p}$, $p \in \mathbb{Z}^+$.
 For $0 < k < p$,
 let $S$ be the $U(1)$-bundle associated with $L^{-k}$, which is a regular Sasaki manifold.
 Let $Z$ be the zero section of $L^{-k}$.
 Suppose that
 $S$ admits a Sasaki-Einstein metric (or $M$ admits a K-E metric). Then there exist
 shrinking and expanding solitons on $L^{-k} - Z$, and they can be pasted together to form
  an eternal solution of the K\"ahler-Ricci flow on
 $(L^{-k} - Z) \times (-\infty, \infty)$. The shrinking soliton for $t \in (-\infty, 0)$ extends smoothly to the zero section $Z$.
 \end{theorem}

Our construction generalizes the work of H.-D. Cao \cite{Cao}, Feldman-Ilmanen-Knopf \cite{Feldman-Ilmanen-Knopf}  which are rotationally symmetric. B. Yang \cite{Yang} and Dancer-Wang \cite{Dancer-Wang} constructed cohomogeneity-one examples. All these examples have very large symmetry group. In contrast, our examples do not carry any continuous symmetry in general.

\section{The Lagrangian mean curvature flow}

An $n$-dimensional submanifold $L\subset \C^n$ is Lagrangian if $\omega|_L=0$ where $\omega$ is the standard K\"ahler form of $\C^n$. An $n-1$ dimensional submanifold $\Sigma \subset S^{2n-1}$ is Legendrian if the cone $L$ over $\Sigma$ is Lagrangian. $L$ is minimal in $\C^n$ if and only if  $\Sigma$ is minimal in $S^{2n-1}$. Being Lagrangian is a closed condition preserved by the mean curvature flow \cite{smoczyk} and thus it makes sense to consider the Lagrangian mean curvature flow. This is similar to being K\"ahlerian is a closed condition that is preserved by the Ricci flow.

Consider a  minimal Legendrian $X:\Sigma^{n-1}\rightarrow S^{2n-1}\subset \C^n$ and a curve $\gamma(s):I\rightarrow \C^*$. An Ansatz for Lagrangian submanifold can be taken to be $ \gamma(s)\cdot X: I\times \Sigma^{n-1}\rightarrow \C^n$ where $\cdot$ denotes complex multiplication. The soliton equation for the Lagrangian mean curvature flow  reduces to an ODE for $\gamma(s)$ (Angenent, Ilmanen, Anciaux \cite{Anciaux},\cite{Anciaux-Castro-Romon} ).

An Ansatz which differs from the minimal Legendrian one but still fits into the separation of variable method is the quadric Ansatz.
Here is an example of a soliton of the Lagrangian mean curvature flow which was constructed in  \cite{Lee-Wang1}. Consider a  quadric in $\R^n$ given by
\[\Sigma= \{(x_1,\cdots, x_n)\,|\, \lambda_1 (x_1)^2+\cdots \lambda_n (x_n)^2=1\} \] where
$\lambda_i$ non-zero and $\Sigma_{i=1^n}\lambda_i>0$. Then
\[ L=\{(x_1 e^{i\lambda_1s},\cdots x_n e^{i\lambda_n s} \,|\, (x_1, \cdots, x_n)\in \Sigma, s\in I\}\] is a soliton solution.

 In fact, the trace of the corresponding flow in $\R^n$ is given by the family of hypersurfaces defined by
 \[\Sigma_{i=1}^n \lambda_i x_i^2=(-2t) \sum_{i=1}^n\lambda_i.\] We note that the flow is shrinking for $t<0$ and expanding for $t>0$. The solution for $t\in (-\infty, \infty)$ forms a weak eternal solution of the mean curvature flow  and the singularity at $t=0$ corresponds to a neck pinching with topological change.

The quadric ansatz was used earlier by Joyce \cite{Joyce1, Joyce2} to construct special Lagrangians submanifolds. Joyce-Lee-Tsui in \cite{Joyce-Lee-Tsui} constructed translating solitons of the Lagrangian mean curvature flow of the form
\[ L=\{(x_1 w_1(s),\cdots x_n w_n(s) \,|\, (x_1\cdots, x_n)\in \Sigma, s\in I\},\] where$w_1(s),\cdots w_n(s)$ satisfy an ODE system. This is a Lagrangian version of ``Grim Reaper" for curve shortening flows. Further generalizations were obtained in \cite{Nakahara}. There are also uniqueness results for minimal Lagrangians and Lagrangian 
solitons, see \cite{Lotay-Neves}, \cite{IJS}, \cite{Castro-Lerma1, Castro-Lerma2}.

\section{Generalizations}
 In this section, we discuss a generalization of both the minimal Legendrian and the quadric Ansatz. Take an $(n-1)$ dimensional submanifold $\Sigma$ of $\C^n$ that is isotropic and consider the image of 
$\mathfrak{A}(s) \Sigma$ where $\mathfrak{A}(s)$ is a curve in the complex affine group which is the semi-direct product of $GL(n, \C)$ and $\C^n$.

One can ask for what kind of $\Sigma$  there will exist a solution $\mathfrak{A}(s)$ such that $\mathfrak{A}(s)\Sigma$ is a Lagrangian submanifold.
The answer is that $\Sigma$ is an \textit{initial data} in the following sense:

A $k$-dimensional isotropic submanifold $\Sigma$ of $\C^n$ is said to be an \textit{initial data} with respect to $(B, b)$ with $B\in GL(n, \C)$ and ${b}\in \C^n$ if at any point $p\in \Sigma$,
\begin{equation}\label{initial_orthogonal} \langle V, BX_p+b\rangle=0\end{equation} for any $V\in T_p\Sigma$, where $\langle\cdot, \cdot\rangle$ is the Hermitian product on $\C^n$.

 For example,  let $\Lambda$ be an $n\times n$  Hermitian matrix and let $\mathfrak{H}$ be the real hypersurface in $\C^n$ defined by $\{X\,\,| \,\, Re\langle \Lambda X, X\rangle=\text{c}\}$ for a constant $c$. Then \[\gamma(\cdot)=\omega(\Lambda X, \cdot) \] is a contact form.
Let $\Sigma$ be a Legendrian submanifold of $\mathfrak{H}$ with respect to $\gamma$, then $\Sigma$ is an initial data.

\begin{theorem}
Given any initial data $(\Sigma, B, {b})$ and any smooth complex curve $\alpha: I\rightarrow \R$, suppose $(A(s), {a}(s))$ satisfies
 \begin{equation}\begin{split}\label{dot_A} A^*(s) \dot{A}(s) &=\alpha(s) B\\
 A^*(s) \dot{a}(s)&=\alpha(s) b \text{ with }  A(0)=I \text{ and } a(0)=0 \end{split}\end{equation} then $A(s)\Sigma+a(s)$ is an isotropic submanifold.
\end{theorem}

In order to reduce the equation to an ODE, it is necessary that the Lagrangian angle of $\mathfrak{A}(s)\Sigma$ depends only on the parameter $s$.
Suppose $(\Sigma, B, b)$ is an initial data of dimension $n-1$ in $\C^n$. For any complex curve $\alpha(s)$, the Lagrangian angle of $\mathfrak{A}(s)\Sigma=A(s)\Sigma+a(s)$ at $(s,p)$ is given by
\begin{equation}\theta (s, p)=arg \det A(s)+arg \alpha(s)+\theta (0, p).\end{equation}

Therefore,  the Lagrangian angle of $\mathfrak{A}(s)\Sigma$ depends only on the parameter $s$ if $\Sigma$ is an \textit{initial data with constant Lagrangian angle} in the following sense: We say $(\Sigma, B, b)$ is an \textit{initial data with constant Lagrangian angle} if the Lagrangian angle of $T_p\Sigma\oplus \R(BX_p+b)$ is a constant in $p\in \Sigma$.

\begin{example}
Suppose $\Sigma$ is an $(n-1)$ dimensional Legendrian submanifold of $\mathfrak{H}=\{X\,\,|\,\, \langle \Lambda X, X\rangle=\text{c}\}\subset \C^n$. The $\Sigma$ is an initial data with constant Lagrangian angle if
\begin{equation}\label{cont_hyper}\omega(-H+\frac{\Lambda^2 X}{|\Lambda X|^2}, \cdot)=0\end{equation} on $\Sigma$ where $H$ is the mean curvature vector of $\Sigma$ in $\mathfrak{H}$.
\end{example} Both minimal Legendrians of $S^{2n-1}$ ($\Lambda=I$) and real quadrics in $\R^n$ ($\mathfrak{H}\cap \R^n$ for $\Lambda$ diagonal) satisfy this condition.

 For an initial data with constant Lagrangian angle $\Sigma$ , a mean curvature equation (minimal Lagrangian equation or
 Lagrangian soliton equation) on the Lagrangian submanifold $A(s)\Sigma$ reduces to an ODE system on $A(s)$ that can be solved.

\begin{theorem}
Suppose $(\Sigma, B, b)$ is an initial data with constant Lagrangian angle. Let $(A(s), a(s))$ be a solution of  \begin{equation}\begin{split}\label{ode} A^*(s) \dot{A}(s) &=\overline{\det A(s)} B\\
 A^*(s) \dot{a}(s)&=\overline{\det A (s)} b \text{ with }  A(0)=I \text{ and } a(0)=0 \end{split}\end{equation} then $A(s)\Sigma+a(s)$ is a special Lagrangian submanifold.
\end{theorem}

We also consider the multiple developing case, i.e. \[A(s_1, \cdots s_k)=A_1(s_1)\cdots A_k(s_k)\] on an $(n-k)$ dimensional initial data and the twisted product of two initial data.
 \begin{theorem} For $p>1$ and $q>1$.
 Suppose $\Sigma_1$ is a $p-1$ dimensional isotropic submanifold of $\C^p$ that is an initial data of constant Lagrangian angle with respect to $B_1$ and $\Sigma_2$ is a $q-1$ dimensional isotropic
submanifold of $\C^q$ that is an initial data of constant Lagrangian angle with respect to $B_2$, then $\Sigma=\Sigma_1\times \Sigma_2\subset \C^{p+q}$ is an initial data of constant Lagrangian angle in $\C^{p+q}$ with respect to any  $c_{11} \tilde{B}_1+c_{12} \tilde{B}_2, c_{21}\tilde{B}_1+c_{22}\tilde{B}_2$ where $\tilde{B}_1=\left[\begin{matrix} B_1&0_{p\times q}\\0_{q\times p} &0_{q\times q}\end{matrix}\right]$,$\tilde{B}_2=\left[\begin{matrix} 0_{p\times p} &0_{p\times q}\\0_{q\times p} &B_2\end{matrix}\right]$, and $C=\left[\begin{matrix} c_{11}&c_{12}\\c_{21}&c_{22}\end{matrix}\right]\in GL(2, \R)$.
\end{theorem}

The twisted product of Legendrians by Castro-Li-Urbano \cite{Castro-Li-Urbano} (see also Haskins-Kapouleas \cite{Haskins-Kapouleas}) is such an example. New examples can be constructed by taking $\Sigma_1$ to be a real quadric in $\C^p$ and $\Sigma_2$ to be a minimal Legendrian in $S^{2q-1}\subset \C^q$.

 Further questions to be studied include:
 
 1. The construction in the last section only addresses the local solvability of the soliton equation. It would be interesting to see whether these solutions
 are complete.
 
 2. Are there more general initial data sets as ``minimal Legendrians" of general real codimensional 1 hyperquadrics in $\C^n$?
 
 3. What is the relation to complex affine geometry? There should be a rigidity theorem for those solitons whose level sets of Lagrangian angles are real hyperquadrics in a Lagrangian subspace? Note that the level sets are congruent under actions of the complex affine group.
 
 4. Is there a similar construction of ``initial data" in the K\"ahler-Ricci flow case? Does it lead to a further generalization of $\eta$-Sasaki-Einstein manifolds?

\end{document}